%

\magnification = 1200

\raggedbottom  
\input amssym.def
\input amssym.tex


\font\bigtenrm = cmr10 scaled \magstep1

\def\rest{\upharpoonright}  

\def\eop{$\bigstar$}  

\def\1st{1$^{\hbox{st}}$}
\def\2nd{2$^{\hbox{nd}}$}
\def\MA{$MA_{ma}(\omega_1)$ }

\def\SS{{\cal S}}
\def\TT{{\cal T}}
\def\MM{{\cal M}}

\def\BB{{\cal B}}

\def\NN{{\cal N}}
\def\CC{{\cal C}}
\def\II{{\cal I}}
\def\FF{{\cal F}}
\def\PPP{{\cal P}}

\def\cc{{\frak c}}

\def\PP{{\Bbb P}}
\def\RR{{\Bbb R}}

\def\sqle{\sqsubseteq}



\def\today{\ifcase\month\or
   January\or February \or March \or April \or May \or June \or
   July \or August \or September \or October \or November \or December \fi
   \space \number\day, \number\year}

\centerline{{\bigtenrm Properties of the Class of Measure Separable
Compact Spaces}
}

\centerline{Submitted to Fundamenta Mathematicae 10 August, 1994}

\centerline{Mirna D\v zamonja and 
Kenneth Kunen\footnote{${}^1$}{\sevenrm Authors
supported by NSF Grant DMS-9100665. First author also partially supported by
the
Basic Research Foundation Grant number 0327398 administered by the
Israel Academy of Sciences and Humanities.}}

\centerline {\sevenrm Hebrew University of Jerusalem\ \ \ \& 
\ $\,$ University of Wisconsin\ \ \ \ \ \ \ \ \ \ \ }
\centerline {\sevenrm \ \ 91904 Givat Ram$\,$\ \ \ \ \ \ \ \ $\,$Madison, WI  53706}
\centerline{\sevenrm \ Israel$\,$\ \ \ \ \ \ \ \ \ U.S.A.}
\centerline {\sevenrm $\,$ dzamonja@sunset.huji.ac.il\ \ $\,$  \& \ \ 
kunen@cs.wisc.edu\ \ \  \ \ \ \ \ \ \ \ \ \ }

\centerline{August 10th, 1994}

\bigskip

\vskip 2 true cm

\centerline{\bf Abstract}

\smallskip
{\narrower\narrower
We investigate properties of the class of compact spaces on which
every regular Borel measure is separable. This class will be referred
to as $MS$.

We discuss some closure properties of $MS$, and show that some simply defined
compact spaces, such as compact ordered spaces or compact scattered spaces,
are in $MS$. Most of the basic theory for regular measures
is true just in $ZFC$. On the other
hand, the existence of a compact ordered scattered space which carries a
non-separable (non-regular) Borel measure
is equivalent to the existence of a real-valued
measurable cardinal $\leq\frak c$.

We show that
not being in $MS$ is preserved by all forcing
extensions which do not collapse $\omega_1$, while
being in $MS$ can be destroyed even by a $ccc$ forcing.

\smallskip

\smallskip
}

\vskip 2 true cm

{\bf \S0. Introduction.} As we learn in a beginning measure theory course,
every  Borel measure on a compact metric space is separable. 
It is natural to ask to what extent this generalizes to other compact spaces. 
It is also true that every Borel measure 
on a compact metric space is regular. 
In this paper, we study the class,  $MS$, of
compacta, $X$, with the property that every
{\it regular \/} measure on $X$ is separable. This contains
some simple spaces (such as compact ordered spaces and compact scattered
spaces), and has some interesting closure properties.
One might also study the class of compacta $X$ such that every 
{\it Borel\/} measure on $X$ is separable,
but the theory here is very sensitive to the axioms of set theory;
for example, the existence of an ordered scattered compactum with
a non-separable Borel measure is independent of $ZFC$
(see Theorem 1.1).  There is also extensive literature about compacta in
which all Borel measures are regular [5].

For the class  $MS$, defined using {\it regular\/} measures, there are still
some independence results, but most of the basic theory goes through in $ZFC$.

First, some definitions:

All spaces considered here are Hausdorff.

We shall consider primarily finite Borel measures on compact spaces.

If $\mu$ is a Borel measure on $X$, the measure
algebra of $(X,\mu)$ is the Boolean algebra of all Borel sets modulo
$\mu$-null sets.  If $\mu$ is finite, then
such a measure algebra is also a metric space, with
the distance between two sets being the measure of their symmetric
difference.  Then, we say that $\mu$ is {\it separable\/} iff
this metric space is separable as a topological space.

A Borel measure $\mu$ on $X$ is {\it
Radon\/} iff the measure of compact sets is finite and 
the measure of each Borel set is the supremum of the measures of its
compact subsets.  If $X$ is compact, this implies
that the measure of each Borel set is also the infimum of the measures of its
open supersets.
Note that for compact spaces, the Radon measures are simply the regular
Borel measures.

The {\it Baire\/} sets are the sets in the least $\sigma$-algebra
containing the open $F_\sigma$ sets.  If $X$ is compact and
$\mu$ is a finite measure defined on the Baire sets, then
$\mu$ extends uniquely to a Radon measure (see [8], Theorem 54D),
and every Borel set is equal to a Baire set modulo a null set.

\medskip
{\bf Definition.}  $MS$ is the class of all compact spaces
$X$ such that every Radon measure on $X$ is separable.
\medskip

Observe, by the above remarks, that if $X$ is compact,
then $X$ is in $MS$ iff every finite Baire measure on X is separable.
We shall primarily be concerned with properties of $MS$,
but we shall occasionally (see Theorem 1.1) remark on finite
non-regular Borel measures, in which case non-separability could
arise from a large number of non-Baire Borel sets.

If not specified otherwise, we give  $[0,1]$ and  $2=\{0,1\}$
their usual probability measures, and then 
$[0,1]^J$ and $2^J$ have the usual product measures.  These
measures defined in the usual way would be defined on the
Baire sets, but they then extend to Radon measures.
These product measures
are in fact {\it completion regular} --
that is, for every Borel set $E$, there are Baire $A,B$
such that $A \subseteq E \subseteq B$ and $B \backslash A$ is
a null set -- but we do not need this fact here.

Note that the measure algebras of $2^J$ and $[0,1]^J$ are isomorphic
whenever $J$ is infinite, and they are separable iff $J$
is countable.  So, for uncountable $J$, $2^J$ and $[0,1]^J$
are simple examples of compact spaces which are not in $MS$.

If $\mu$ is a Borel measure on $X$, and $E$ is a Borel set,
then $\mu\rest E$ is the Borel measure on $E$ defined in the obvious
way: $(\mu\rest E) ( B) = \mu(B)$ for Borel $B \subseteq E$.
We say that $\mu$ is {\it nowhere separable\/} iff 
$\mu \rest E$ is non-separable for each Borel set $E$ of
positive measure.

Our basic notions never assume that non-empty open sets
have positive measure, but it is sometimes useful to reduce
to this situation.  If $\mu$ is a Radon measure on the compact
space $X$, let $U$ be the union of all open null sets.
By regularity of the measure, $U$ is also a null set, and is
hence the largest null set.  We call $K = X \backslash U$ the
{\it support\/} of $\mu$.  Note that $\mu(K) = \mu(X)$, 
and every relatively open non-empty subset of $K$ has positive measure.

The following lemma is sometimes useful to reduce
the study of non-separable measures to nowhere separable measures:

\medskip

{\bf Lemma 0.0.}  If $X$ is compact and $\mu$ is a non-separable
Radon measure on $X$, then there is a closed $K \subseteq X$
such that $\mu(K) > 0$, $\mu\rest K$  is nowhere separable,
and every relatively open non-empty subset of $K$ has positive measure.

{\bf Proof.}  By  Maharam's Theorem [13],
there is a Borel $E \subseteq X$
such that $\mu(E) > 0$ and $\mu\rest E$  is nowhere separable.
We then apply regularity of $\mu$ to choose $C \subseteq E$
of positive measure,
and let $K$ be the support of $\mu \rest C$.  \eop

\medskip

In \S1, we consider some classes of topological spaces which
are subclasses of
$MS$, and
in \S2, we discuss various closure properties of $MS$.

In \S\S 3,4, we look at the behavior
of $MS$ in transitive models of set theory.
It is easy to see that the property of {\it not\/} being in $MS$
is preserved under any
forcing extension which does not collapse $\omega_1$.
In \S4, we
show that being in $MS$ need not even be preserved by $ccc$ forcing;
assuming the existence of a Suslin tree $T$, 
we construct an $X \in MS$ such 
that forcing with $T$ adds a non-separable Radon
measure on $X$ in the generic extension.
Of course, since the notion of ``compact space'' is not absolute
for models of set theory,
some care must be taken to say precisely what
is meant by looking at the same compact space in two different
models; this is handled in \S3, and in a somewhat different way
by Bandlow [1].

We do not know if there is any simple way of expressing 
``$X \in MS$'' without mentioning measures.
By the results of \S\S 1,2, there are some simple sufficient 
conditions for a compact space $X$ to be in $MS$; for example,
it is sufficient that $X$ be a subspace of a countable product of
ordered spaces and scattered spaces.  By the result of \S4,
any condition of this form, which is preserved in the passage to
a larger model of set theory, cannot be a necessary condition 
as well (or, at least, cannot be proved to be necessary in ZFC).

\bigskip

{\bf\S1. Subclasses of $MS$.}  We begin by pointing out
some simple sufficient conditions for a compact space to be in $MS$.

First, recall some definitions.
A topological space is {\it ccc} iff there are no uncountable
disjoint families of open subsets of the space.
If $\mu$ is a Radon measure on a compact space, $X$, then
$X$ need not be ccc, but the support of $\mu$ is ccc.
A space $X$ is a {\it LOTS\/} (linearly ordered
topological space) if its topology is the
order topology induced by some total order on $X$.

\medskip

{\bf Theorem 1.0.}  $MS$ contains every compact $X$ such that $X$ satisfies
one of the following.
{\parindent=30pt
 \item{1.} $ X$ is second countable (= metric).

 \item{2.}$ X$ is scattered.

 \item{3.} Every $ccc$ subspace of $X$ is second countable.

 \item{4.} $X$ is Eberlein compact.

 \item{5.} $X$ is a LOTS.
}

\medskip

{\bf Proof.} Suppose that $X$ is  compact and $\mu$ is a Radon
measure on $X$.

For (1), 
fix a countable basis $\BB$ for $X$, which is
closed under finite unions, and note that $\BB$ is dense in the
measure algebra of $(X,\mu)$.
For (2) and (3), if $\mu$ were non-separable, then the $K$
provided by Lemma 0.0 would yield an immediate contradiction.
Now, (4) follows because, by Rosenthal [14], 
every $ccc$ Eberlein compact is second countable.

For (5), assume that $X$ is a compact LOTS and that $\mu$ is
non-separable.  By Lemma 0.0 and the fact that every
closed subspace of $X$ is a LOTS, we may assume, without
loss of generality, that $\mu$ is nowhere
separable on $X$; in particular, every point in $X$ is
a null set.  We may also assume that $\mu(X)=1$.
Let $a$ be the first element of $X$ and $b$ the last element of $X$.
Define $f : X \to [0,1]$ by: 
$f(x)=\mu([a,x])$.  Then $f$ is continuous (since points are null sets),
$f(a) = 0$, and $f(b) = 1$. Let $\lambda=\mu f^{-1}$ be the
induced Borel measure on $[0,1]$.  Then $\lambda$ is regular and separable.
Let $\Sigma$ be the family of all Borel subset $B$ of $X$
such that there is a Borel
subset $E$ of $[0,1]$ with $\mu(B\Delta f^{-1}(E))=0$.
To conclude that $\mu$ is separable (and hence a contradiction),
it is sufficient to show that $\Sigma$ in fact contains all Borel sets, since
then the measure algebras of $(X,\mu)$ and $([0,1], \lambda)$
will be isometric.
This will follow if we can show that $\Sigma$ contains
all Baire sets. 
Since $\Sigma$ is
a $\sigma$-algebra and every Baire set is in 
the $\sigma$-algebra generated by intervals, it
is sufficient to show that $\Sigma$ contains all intervals.
Since $\Sigma$ certainly contains all singletons (take $E = \emptyset$),
it  is sufficient to show that each $[a,x] \in \Sigma$.
Fix $x$, and let $s=f(x)$, and $E = f([a,x]) = [0,s]$;
then $f^{-1}(E) = [a,z]$ for some $z \geq x$ with $f(z) = s$.
$[a,x] \subseteq [a,z]$, and $\mu([a,x]) = f(x) = f(z) = \mu([a,z])$,
so $\mu([a,x] \Delta [a,z]) = 0$. \eop

\medskip

The proof of (5) would have been a little nicer if we could
have said that $f$ were 1--1, since that would have implied
that $X$ is second countable.  But we cannot say this.
Even if all non-empty open subsets of $X$ have positive measure,
there could be points $x < z$ with no points between them, 
in which case $f(x) = f(z)$.  For a specific example, take
$X$ to be the double arrow space, which is not second
countable but which is the support of a Radon measure.

Regarding (4), the statement that all
Corson compacta are in $MS$ is independent of $ZFC$.
See Kunen and van Mill [12] and \S2 for further discussion.

The proofs of (2) and (5) involve passing to the  support
of the measure, by Lemma 0.0, which is justified by regularity
of the measure.
If we drop regularity, $X$ can be both scattered {\it and\/} a LOTS
and still have a non-separable Borel measure:

\medskip

{\bf Theorem 1.1.} There is a compact scattered  LOTS which has
a non-separable finite Borel measure iff there is
a real-valued measurable cardinal $\le \cc$.

{\bf Proof.} If $\kappa$ is real-valued measurable,
let $\mu$ be a real-valued measure on $\kappa$ such that
the set of limit ordinals is a null set; then every subset 
of $\kappa$ is equal to a Borel (in fact, open) set modulo
a null set.  This measure is
non-separable by the Gitik-Shelah Theorem [4,6,7].
So,  $\mu$ on $\kappa+1$ yield
an example of an ordered scattered continuum having a
non-separable Borel measure.

Now we show that,
if there are no real-valued measurable cardinals $\le \cc$,
and $\mu$ is a finite Borel measure on a compact
scattered LOTS $X$, then $\mu$ is completely atomic.

We do not lose generality if we
assume that $\mu$ is atomless on $X$, 
there are no real-valued measurable cardinals $\le \cc$,
and $\mu(X) = 1$.  We derive a contradiction.

{\it Remark\/}:  If $S \subseteq X$ has the property that every
subset of $S$ is Borel, then $\mu(S) = 0$
(by no real-valued measurable cardinals).  More generally,
call $(S,f,\theta)$ a {\it dangerous triple\/} iff
$S$ is a Borel subset of $X$,  $\mu(S) > 0$,  $\theta$ is a cardinal, 
and $f : S \to \theta$ has the property that $f^{-1}(Z)$ is Borel
for each $Z \subseteq \theta$ and $\mu(f^{-1}(\{z\})) = 0$
for all $z \in \theta$.  Then the induced measure, $\mu f^{-1}$, is
a non-trivial measure defined on all subsets of $\theta$,
and must then be completely atomic (again,
by no real-valued measurable cardinals).  This
is not immediately a contradiction (unless there are no two-valued
measurable cardinals either).  But, since $\mu$ is atomless,
there must be a Borel $Y \subseteq S$ which is not equal
to any $f^{-1}(Z)$  (for $Z \subseteq \theta$) modulo a null set.
We shall use this remark later.

Let $X^{(\alpha)}$ be the $\alpha^{\rm th}$ derived subset of $X$.
If $x \in X$, let $rank(x)$ be the least $\alpha$ such
that $x \notin X^{(\alpha+1)}$.
If $C$ is a non-empty closed subset of $X$, let $rank(C)$
be the least $\alpha$ such
that $C \cap X^{(\alpha+1)} = \emptyset$.
Note that if $\alpha = rank(C)$, then
$C \cap X^\alpha$ is finite and non-empty.

Let $\CC$ be the set of all closed $C \subset X$ such that
$\mu(C) = 0$ and $C$ contains the first and last elements of $X$.
If $C  \in \CC$, let $\II(C)$ be the set of all non-empty maximal
intervals of $X \backslash C$.  
If $x \in X \backslash C$, let $comp(x,C)$ be the (unique) 
$I \in \II(C)$ such that $x \in I$.  Note that if $C,D \in \CC$,
$C \subseteq D$, and $x \in X \backslash D$, then
$comp(x,D) \subseteq comp(x,C)$.

If $C,D \in \CC$, say $C \ll D$ iff $C \subset D$ and for all
$x \in X \backslash D$, $rank(cl(comp(x,D))) < rank(cl(comp(x,C)))$
(here, $cl$ denotes topological closure).
Observe that if we get $C_n \in \CC$ for $n \in \omega$ with
each $C_n \ll C_{n+1}$, we will have a contradiction, since
$\bigcup_{n \in \omega}C_n$ will
have measure 0 and equal $X$ (since an $x$ not in the union
would yield a decreasing $\omega$-sequence of ordinals).

Thus, it is sufficient to fix $C \in \CC$ and find a $D \in \CC$
with $C \ll D$.  First, note that if $S \subseteq X \backslash C$
and $S$ contains at most one point from each $I \in \II(C)$,
then every subset of $S$ is Borel, so $\mu(S) = 0$.
So,  $\mu(S) = 0$ whenever $S$ contains at most countably
many points from each $I \in \II(C)$.

By expanding $C$ if necessary, we may assume that for
each $(a,b) \in \II(C)$,  the points of maximal rank in $[a,b]$
are among $\{a,b\}$.

For each $(a,b) \in \II(C)$:  If $b$ is a successor point, let 
$R_0(b)$ be the singleton of its predecessor.  If
$cf(b) = \omega$, let $R_0(b)$  be some increasing $\omega$-sequence in 
$(a,b)$ converging to $b$.  Otherwise, let $R_0(b) = \emptyset$.
Likewise define $L_0(a)$ to be a singleton if $a$ is a predecessor
point, a decreasing $\omega$-sequence if $ci(a) = \omega$,
and empty if $ci(a) > \omega$.

Let $\FF$ be the set of all closed $D \supseteq C$ such that
$D$ is of the form 
$$
C \cup \bigcup\{ R(b) \cup L(a) : (a,b) \in \II(C)\}\ \ \ ,
$$
where for each $(a,b) \in \II(C)$:  
$R(b) = R_0(b)$ if $R_0(b) \ne \emptyset$, and otherwise 
$R(b)$ is a closed cofinal sequence of type $cf(b)$
in $(a,b)$ converging to $b$; and,
$L(a) = L_0(a)$ if $L_0(a) \ne \emptyset$, and otherwise 
$L(a)$ is a closed coinitial sequence of type $ci(a)$
in $(a,b)$ converging to $a$.

Note that $\FF$ is closed under countable intersections,
so we may fix $D \in \FF$ of minimal measure.   Then, 
note that $\mu(D) = 0$.  To see this, consider 
$(S,f,\theta)$, where $S = D \backslash C$, 
$\theta = |\II(C)|$, and $f$ maps $I \cap (D \backslash C)$
to one point in $\theta$ for each $I \in \II(C)$.
If $\mu(D) > 0$, then $(S,f,\theta)$ would be a dangerous triple.
But also, note that if $cf(b) > \omega$, then every Borel set either
contains or is disjoint from  a closed cofinal sequence in $b$.
Using this, and minimality of $\mu(D)$, we see that every
Borel $Y \subseteq S$ is equal
to some $f^{-1}(Z)$  modulo a null set, which is a contradiction.

So, $C \ll D$.  \eop

\bigskip

{\bf \S2. Closure Properties of $MS$.}  In this section,
we consider questions about the closure of $MS$ under
under subspaces, continuous images, continuous pre-images, and products.
We begin with:

\medskip

{\bf Lemma 2.0.} If $X \in MS$, then every closed subspace of $X$ is in MS.

\medskip

Of course, this is trivial, since a measure on a subspace
induces a measure on $X$ in the obvious way.  The same argument
works for continuous images, but requires a little care:

\medskip

{\bf Lemma 2.1.} Suppose that $X \in MS$ and $f$ is a continuous
map from $X$ onto $Y$. Then $Y$ is in $MS$.  

{\bf Proof.}  Suppose $\mu$ were a non-separable Radon measure
on $Y$.  Choose a Radon measure 
$\nu$ on $X$ such that $\mu=\nu f^{-1}$.
The existence of such a $\nu$ follows easily from the Hahn-Banach
Theorem plus the Riesz Representation Theorem; see also
Henry [10], who proved this, plus some more general results.
Now, the measure algebra of $\mu$ embeds into the
measure algebra of $\nu$, so $\nu$ is non-separable,
contradicting $X \in MS$. \eop

\medskip

In particular, if $X$ maps onto $[0,1]^{\omega_1}$,
then $X \notin MS$.  It is a well-known open question of Haydon whether
the converse holds; that is, 
if $X$ is compact and  $ X \notin MS$, must $X$ map onto $ [0,1]^{\omega_1}$?
Many counter-examples are known under $ CH$ or some other axioms of set
theory [2,9,11,12],
but it is unknown whether a ``yes'' answer is consistent,
or follows from $ MA + \neg CH$.

\smallskip

We shall now show that $MS$ is closed under countable products;
it is obviously not closed under uncountable products.
First, consider a product of two spaces:

\medskip

{\bf Lemma 2.2.} If $X,Y \in MS$, then $X\times Y \in MS$.

{\bf Proof.} Let $\lambda$ be a Radon
measure on $X\times Y$.  We show that $\lambda$ is separable.

Let $\mu$ be the Radon measure on $X$ induced from $\lambda$
by projection on the first co-ordinate.  Since $X \in MS$,
there is a countable family $\{D_n:n\in\omega\}$ of
closed subsets of $X$ which is dense in the measure algebra of $(X,\mu)$. 

For each $n$, let $\nu_n$ be the Radon
measure on $Y$ induced from
$\lambda\rest (D_n\times Y)$ by projection on the second co-ordinate.
Since $Y \in M$, for each $n$ there
is a family $\{ E_m^n: m\in\omega\}$ of closed subsets
in $Y$ which is dense in the measure algebra of $(Y,\nu_n)$.

Then the family of the finite unions of the sets of the form
$D_n\times E_m^n$ is dense in the measure
algebra of $(X\times Y, \lambda)$.\eop

\medskip

{\bf Theorem 2.3.} $MS$ is closed under countable products.

{\bf Proof.} Suppose that $X_n(n\in\omega)$ are in $MS$ and $\mu$ is a
Radon measure on $X=\Pi_{n\in\omega} X_n$.

For every $n$, let $\pi_n$
denote the natural projection from $X$ onto $Y^n=\Pi_{k\le n}X_k$.
Then $\mu_n=\mu\, \pi_n^{-1}$ is a Radon measure on $Y^n$, and therefore
separable, by the previous Lemma (plus induction).
For each $n$, fix a countable family
${\cal D}_n$ which is dense in the measure algebra of $(Y^n,\mu_n)$.
Then ${\cal D}=\bigcup_{n\in\omega} \{ \pi_n^{-1} (D):D\in {\cal D}_n\}$
is dense in the measure algebra of $(X,\mu)$.\eop

\medskip

By the same argument:

\medskip

{\bf Lemma 2.4} $MS$ is closed under inverse limits of countable length.

\medskip

Since $MS$ is closed under countable products and not
closed under uncountable products, it is reasonable to
consider now $\Sigma$-products, a notion between countable
and uncountable products.
Let  $X_\alpha(\alpha\in\kappa)$ be topological spaces, let
$X$ be the usual Tychonov product of the $X_\alpha$,
and let $a=\langle a_\alpha: \alpha\in\kappa \rangle $ be a point in $X$.
We define $\Sigma(a)$ to be the set of all points
$x$ of $X$ which differ from $a$ on just a countable set of coordinates.
Considered as a subspace of $X$, this set is
called the {\it  $\Sigma$-product\/} of the $X_\alpha$ with
{\it base point\/} $a$.
If $\kappa$ is countable, this is just the Tychonov product.
If $\kappa$ is uncountable, then except in trivial
cases, $\Sigma(a)$ is not compact and is a proper subset of the
Tychonov product.  So, the question we address now is:  if
each $X_\alpha \in MS$, must every compact subspace
of $\Sigma(a)$ be in $MS$?  The answer turns out to independent
of $ZFC$, and in fact equivalent to a weakened version of
Martin's Axiom ($MA$).

Let \MA denote the statement that $MA(\omega_1)$ is true for
measure algebras; that is, whenever $\PP$ is a ccc partial
order which happens to be a measure algebra, then
one can always find a filter meeting $\omega_1$ dense sets.
So, \MA implies $\neg CH$, and \MA  follows from
$MA(\omega_1)$.  But also,
\MA  is true in the random real model, or in any model
with a real-valued measurable cardinal, where most of the combinatorial
consequences of full $MA$ fail (see [4]).
Consequences of
$MA(\omega_1)$ for measure algebras in measure theory are numerous (see [3]),
and some of them really only require \MA.

By Kunen and van Mill [12], \MA is equivalent to the fact that
all Corson compacta are in $MS$.  Recall 
that $X$ is called a {\it Corson compact \/} iff
$X$ is homeomorphic to a compact subspace of a $\Sigma$-product
of copies of $[0,1]$.  So, if \MA fails, there is a 
compact subspace of a $\Sigma$-product of spaces in $MS$
which fails to be in $MS$.  Conversely, we can adapt the
proof in [12] to show:

\medskip

{\bf Theorem 2.5.} Assuming \MA,
if $Y$ is a compact subspace of a
$\Sigma$-product of spaces in $MS$, then $Y \in MS$.

{\bf Proof.}  Suppose that $Y$ is a compact
subspace of the $\Sigma$-product of the $X_\alpha$ ($\alpha \in \kappa$),
with base point $a$, where each $X_\alpha \in MS$.
Assume that $\mu$ is a non-separable Radon measure
on $Y$.  By Lemma 0.0, we may assume that every non-empty
relatively open subset of $Y$ has positive measure.
Let $J = \{\alpha  \in \kappa: \exists y \in Y ( y_\alpha \ne a_\alpha) \}$.
If $J$ is countable, then $Y$ is a homeomorphic to a
closed subspace of the
Tychonov product of the $X_\alpha$ ($\alpha \in J$), so $Y$ would
be in $MS$ by  Theorem 2.3 and Lemma 2.0.  So, we
assume $J$ is uncountable and derive a contradiction.

Choose distinct $\alpha_\xi \in J$ for $\xi < \omega_1$.  For each
$\xi$, let $\pi_\xi: Y \to X_{\alpha_\xi}$ be the natural projection.
For each $\xi$, there is a $y_\xi \in Y$ with
$\pi_\xi(y_\xi) \ne a_{\alpha_\xi}$, and hence there is a
relatively open $U_\xi \subseteq Y$ such that 
$a_{\alpha_\xi} \ne \pi_\xi(\overline U_\xi)$.

Since each $U_\xi$ has positive measure, we can apply \MA
to find an uncountable $L \subseteq J$
such that $\{U_\xi : \xi \in L\}$ has the finite intersection
property. $L$ exists because $MA(\omega_1)$ 
for a ccc partial order implies that the order has $\omega_1$
as a precaliber.  Here the order in question is the measure
algebra of $X$.

By compactness, choose 
$z  \in \bigcap_{\xi \in L} \overline U_\xi$.
Then $z_\xi \ne a_\xi$ for all $\xi \in L$, contradicting
the definition of $\Sigma$-product.  \eop

\medskip

We now consider the situation with continuous preimages of spaces in $MS$.
Suppose $X$ is compact, $f: X \to Y$, and $Y \in MS$.  Obviously, we cannot
conclude $X \in MS$, since $2^{\omega_1}$ maps onto $2^\omega$.
But we might hope to conclude $X \in MS$ if we know also that 
the preimage of each point  is in $MS$.   Unfortunately, this
is false, at least under $CH$, by an example of Kunen [11]:
under CH, there is a closed subset $X$ of $2^{\omega_1}$
such that $X$ supports a non-separable Radon probability measure, yet, the
projection $f :X\to 2^\omega$ satisfies that
for each $y\in 2^\omega$, $f^{-1}\{y\}$ is second countable.

However, there are two special cases where we can conclude
from $f: X \to Y$
that $X \in MS$.  One (Theorem 2.7) is where $Y \in MS$ and
the point preimages
are scattered.  The other (Theorem 2.9) is where the point
preimages are in $MS$ and $Y$ is scattered.
Of course, there is a third special case which we have 
already covered: if $X$ is a product, $Y \times Z$, and $f$
is projection, it is sufficient that the point preimages (i.e. $Z$)
be in $MS$ to conclude $X \in MS$ by Lemma 2.2.

In the proof of Theorem 2.7, we shall use the following general notation. 
Suppose $X$ and $Y$ are compact, $f: X \to Y$, and $\mu$ is a Radon measure
on $X$.  Let $\nu = \mu f^{-1}$ be the induced measure on $Y$.
If $E$ is any Borel subset of $X$, let $\nu_E$ be the measure on
$Y$ defined by $\nu_E(B) = \mu(E \cap f^{-1}(B))$.
Clearly, $0 \le \nu_E \le \nu$.  Let $\delta(E) \in L^1(\nu)$ be the 
Radon-Nikodym derivative of $\nu_E$; so $d \nu_E = \delta(E) d \nu$.
Then $0 \le \delta(E)(x) \le 1$ for ae $x$.  In the following,
$\| \cdot \|$ always denotes the $L^1$ norm on $L^1(\nu)$.

The next lemma shows how to split a closed subset of $X$
into two independent pieces.

\medskip

{\bf Lemma 2.6.}  Suppose that
$X$ and $Y$ are compact and $f: X \to Y$.
Suppose that $\mu$ is a nowhere separable Radon measure on $X$, but
$\nu = \mu f^{-1}$ is a separable measure on $Y$.
Let $H \subseteq X$ be closed, and fix $\epsilon >0$.  Then
there are disjoint closed $K_0, K_1 \subseteq H$ such that
for $i = 0,1$,\ \  $\delta(K_i) \le {1 \over 2}\delta(H)$ and
$\| {1 \over 2}\delta(H) - \delta(K_i) \| \le \epsilon$.

{\bf Proof.}  Let $\MM$ be the measure algebra of $X,\mu$.
Let $\NN$ be the sub $\sigma$-algebra of $\MM$ generated by
$H$ and all $f^{-1}(B)$, where $B$ is a Borel subset of $Y$.
Since $\MM$ is nowhere separable while $\NN$ is separable, 
Maharam's Theorem implies that there are complementary Borel
sets $E_0, E_1 \subseteq X$ such that
$\mu(E_0 \cap A) = \mu(E_1 \cap A) = {1 \over 2}\mu(A)$ for
all $A \in \NN$.  In particular, whenever $B \subseteq Y$ is 
Borel, and $i$ is 0 or 1, 
$\mu(E_i \cap H \cap f^{-1}(B)) = {1 \over 2}\mu( H \cap f^{-1}(B) )$.
Thus, $\delta(E_i \cap H) = {1 \over 2}\delta(H)$.

Now, for $i = 0,1$, let $K_i^n$ for $n \in \omega$ be
an increasing sequence of closed subsets of $E_i$, such
that $\mu(K_i^n) \nearrow  \mu(E_i)$.  Then
$\delta(K_i^n) \to \delta(E_i)$ in $L^1(\nu)$, so, for $n$ sufficiently
large, setting $K_i = K_i^n$ will satisfy the Lemma.  \eop

\medskip

{\bf Theorem 2.7.} Suppose that $X$ is compact, $f : X \to Y$,
$Y \in MS$, and $f^{-1}\{y\} $ is scattered for all $y \in Y$.
Then $X \in MS$.

\medskip

{\bf Proof.}  Suppose $X \notin MS$.  We shall find a $y \in Y$
such that $f^{-1}\{y\} $ is not scattered.  Let $\mu$ be a non-separable
Radon measure on $X$.  We may assume that $\mu$ is nowhere separable,
since otherwise we may simply replace $X$ by a closed subset of $X$ 
on which $\mu$ is nowhere separable.

We shall find closed subsets of $X$, $H_s$, for $s \in 2^{< \omega}$,
such that they form a tree:

\smallskip
\item{\bf (1)} $H_{()} = X$.  For each $s$, $H_{s0}$ and $H_{s1}$ are disjoint
non-empty closed subsets of $H_s$.
\smallskip

\noindent
Note, now, that if $ y \in \bigcap\{f(H_s) : s \in 2^{< \omega}\} $,
then $f^{-1}\{y\}$ has a closed subset which maps onto $2^\omega$,
so $f^{-1}\{y\}$ is not scattered.  To ensure that there is such a $y$, we
assume also

\smallskip
\item{\bf (2)} For each $n \in \omega$, there is a closed $L_n \subseteq Y$
such that $f(H_s) = L_n$ for all $s \in 2^n$.
\smallskip

\noindent
Then the $L_n$ will form a decreasing sequence of closed sets,
so, by compactness, we may simply choose $y \in \bigcap_{n \in \omega}L_n$.
So, we are done if we can actually construct the $H_s$ and $L_n$
to satisfy (1,2).    To aid in the inductive construction, we assume
also:

\smallskip
\item{\bf (3)} $\nu(L_n) > 0$ for all $n$.
\item{\bf (4)} $\delta(H_s) \ge 2^{-2n}$ ae on $L_n$, for each $s \in 2^n$.
\smallskip

\noindent
Here, $\nu$ and $\delta(H)$ are as defined above.  Since items
(1-4) are trivial for $n =  0$, we are done if we can show how,
given $L_n$ and the $H_s$ for $s \in 2^n$, we can construct
$L_{n+1}$ and each $H_{s0}, H_{s1}$.  First, apply Lemma 2.6
and choose, for each $s$, disjoint closed
$K_{s0}, K_{s1} \subseteq H_s$ such that
for $i = 0,1$,
$\delta(K_{si}) \le {1 \over 2}\delta(H_s)$ and
$\| {1 \over 2}\delta(H_s) - \delta(K_{si}) \| \le 2^{-3n-4}\nu(L_n)$.
Let 
$$
A_{si} = \{y \in L_n: \delta(K_{si})(y) \le {1 \over 4} \cdot
2^{-2n}\} \ \ \ .
$$
Since ${1 \over 2} \delta(H_s) \ge {1 \over 2} 2^{-2n}$
on $L_n$, 
$$
2^{-3n-4}\nu(L_n) \ge \| {1 \over 2}\delta(H_s) - \delta(K_{si}) \| \ge
\nu(A_{si}) \cdot {1\over 4} 2^{-2n} \ \ \ ,
$$
so $\nu(A_{si}) \le 2^{-n-2} \nu(L_n)$.
Let $B = \bigcup\{A_{si} : s \in 2^n , i = 0,1\}$.
Then $\nu(B) \le {1 \over  2} \nu(L_n)$, so 
$\nu(L_n \backslash B) > 0$.  For all $y \in 
L_n \backslash B$,\ \  $\delta(K_{si})(y) \ge 2^{-2(n+1)}$ for each $s,i$.
In particular, then, $\nu(L_n \backslash B \backslash f(K_{si})) = 0$. 
So, we may
choose a closed $L_{n+1} \subseteq (L_n \backslash B)$ such that 
$\nu(L_{n+1}) > 0$ and $L_{n+1} \subseteq f(K_{si})$ for each $s,i$.
Finally, let $H_{si} = K_{si} \cap f^{-1}(L_{n+1})$.  Note that
$\delta(H_{si}) = \delta(K_{si}) \ge 2^{-2(n+1)}$ on $L_{n+1}$.  \eop

\medskip

Now, before turning to the case that $Y$ is scattered,
let us pursue the following idea.
If $X \notin MS$, $X$ could still have a clopen subset in $MS$; for example,
$X$ could be the disjoint sum of $2^{\omega_1}$ and $2^{\omega}$.
However, if one deletes all the open subsets of $X$ which are in $MS$, one gets
a ``kernel'' which is everywhere non-$MS$ by Theorem 2.8.(d) below.

Given a compact $X$, define
$$
ker(X) = X \backslash \bigcup\{U \subseteq X :
U \hbox{\ is open and\ } cl(U) \in MS\} \ \ \ .
$$

\medskip

{\bf Theorem 2.8.} If $X$ is any compact space:

\item{a.} $ker(X)$ is a closed subset of $X$.
\item{b.} If $Y$ is any closed subset of $X$,
then $ker(Y) \subseteq ker(X)$.
\item{c.} $X \in MS$ iff $ker(X) = \emptyset$.
\item{d.} $ker(ker(X)) = ker(X)$.

{\bf Proof.} (a) is obvious.  (b) follows from Theorem 2.0.
For (c), if $ker(X) = \emptyset$, then by compactness, $X$ is a finite
union of closed sets in $MS$, which clearly implies that $X \in MS$.

If (d) fails, fix $p \in ker(X) \backslash  ker(ker(X))$.
Applying the definition of $ker$ to $ker(X)$,
$p$ has a neighborhood $U$ in $X$ such that
$cl(U \cap ker(X)) \in MS$; let $V$ be a neighborhood of $p$ in $X$
such that $cl(V) \subseteq U$; then (by Theorem 2.0),
$cl(V) \cap ker(X) \in MS$.  Since $cl(V) \notin MS$,
let $\mu$ be a non-separable Radon measure on $cl(V)$.
Applying Lemma 0.0, 
let $K$ be a closed subset of $cl(V)$ such that $\mu(K) > 0$,
$\mu \rest K$ is nowhere separable, and
every relatively non-empty relatively open subset
of $K$ has positive measure.
Then, $K = ker(K)$, and, applying (b), $ker(K) \subseteq ker(X)$,
so $K \subseteq cl(V) \cap ker(X)$, contradicting
$cl(V) \cap ker(X) \in MS$.  \eop

\medskip

{\bf Theorem 2.9.} Suppose $X$ and $Y$ are compact,
$Y$ is scattered, $f: X \to Y$, and the
preimages of all points in $Y$ are in $MS$. Then $X$ is in $MS$.

{\bf Proof.} If $X \notin MS$, $ker(X) \ne \emptyset$, so let
$y \in f(ker(X))$ be an isolated point in $f(ker(X))$.
Then $f^{-1}(y) \cap ker(X)$ is a clopen subset of $ker(X)$,
so $f^{-1}(y) \cap ker(X) \notin MS$ by Theorem 2.8(d),
so $f^{-1}(y)  \notin MS$ by Theorem 2.0.  \eop

\medskip

{\bf Corollary 2.10.} Suppose $S$ is a direct
sum of compact spaces $X_\alpha$, for $\alpha \in \kappa$
(so $S$ is locally compact). Suppose that each $X_\alpha \in MS$.
Then any compactification
of $S$ with remainder in $ MS$ is in $MS$ (in particular, the 1-point
compactification).

{\bf Proof.} Apply Theorem 2.9 with $Y$ being
the 1-point compactification of a discrete $\kappa$,
and $f$ taking each $X_\alpha$ to $\alpha$ and the remainder
to the point at infinity. \eop

\bigskip

{\bf \S3. Compact Spaces in Models of Set Theory.}
In forcing, we frequently discuss the preservation of a property
(such as $MS$) as we pass between two models of set theory.
Suppose that $M \subseteq N$ are two transitive models of $ZFC$,
with $X$ a topological space in $M$.  If $M$ thinks that
$X$ has some property, we may ask whether $N$ also thinks
that the {\it same space\/} $X$ has that property. But,
since being a space is not a first-order notion, we must
be more precise about what ``same space'' means.  There are actually
two possible meanings to this, only one of which makes sense in the
case of $MS$.

The first meaning is the most common one, and is frequently 
used without comment.
Formally, a space is a pair, $\langle X, \TT \rangle$,
where $X$ is a set and $\TT$ is a topology on $X$.  
If $\langle X, \TT \rangle \in M \subseteq N$, and 
the statement ``$\TT$ is a topology on $X$'' is true in $M$,
then this statement will not in general be true in $N$, but it
will be true in $N$ that $\TT$ is a basis for a topology, $\TT'$, on $X$.
In the future, we shall often suppress explicit mention of
$\TT$ and $\TT'$, and simply say something like:  ``Let $X$ be a space
in $M$, and now consider the same $X$ in $N$''.

However, in dealing with properties of compact spaces, such as $MS$, it
is really the second meaning which is required.
If $X$ is a compact space in $M$ (i.e., the statement, 
``$\langle X, \TT \rangle$ is compact'' is true relativized to $M$),
then the same $X$ in $N$ is not necessarily compact. 
For example, if  $X$ is $[0,1]^M$ (the unit interval of $M$),
and $N$ has new reals which are not in $M$, then the same $X$ in $N$ is not
compact from the point of view of $N$; more generally, if $N$
has new reals, then it is only the scattered compact spaces of
$M$ which remain compact in $N$.
If $X$ is a compact space in $M$, we shall define  a compact space in $N$,
which we shall call
$\Phi_{M,N}(X)$, or just $\Phi(X)$ when $M,N$ are clear from context.
Informally, $\Phi(X)$ will be the compact space in $N$ which 
``corresponds'' to $X$.
In some simple cases, $\Phi(X)$ is the ``obvious thing''.  For example,
if  $X$ is the unit interval of $M$,
then $\Phi(X)$ is the unit interval of $N$;
if  $X$ is the $n$-sphere in $M$,
then $\Phi(X)$  is the $n$-sphere in $N$;
if $X$ is a Stone space of a Boolean algebra $B \in M$,
then $\Phi(X)$ is the Stone space of the same $B$ as computed
within $N$. 
But, we must be careful to check that this $\Phi(X)$
is computed for {\it every\/} compact $X$ in some natural way.
Here, ``natural'' can be expressed formally in terms of categories.
Let $ CT$ be the category of compact $T_2$ spaces and continuous
maps.  If $M$ is a transitive model of $ZFC$, $CT^M$ is just the relativized
$CT$, computed within $M$.
Then, $\Phi_{M,N}$ will be a functor from $CT^M$ to $CT^N$.

$\Phi(X)$ will in fact be computed in $N$ as
some compactification of $X$,
so we pause to make some remarks on compactifications.  Here,
we just work in $ZFC$, forgetting temporarily about models.

\smallskip

Let  $C(X)$ denote the family of all bounded continuous real-valued
functions on $X$.  This is a Banach space, and we let $\|f\|$ denote
the usual $\sup$ norm.  Also, $C(X)$ is a Banach algebra under
pointwise product.  If $\SS$ is any non-empty subset of $C(X)$, let
$e_{\SS}$, or just $e$, denote the usual evaluation map
from $X$ into the cube,  $\prod\{[-\|f\|, +\|f\| ] : f \in \SS\}$; that is,
$(e(x))(f) = f(x)$.  Let $\beta(X,\SS)$ be the closure of $e(X)$
in this cube.  It is always true that $e$ is continuous.
In some cases (for example, if $\SS$ separates points and closed sets),
$e$ will be a homeomorphic  embedding of $X$, in which case
$\beta(X,\SS)$ is a compactification of $X$.  If $\SS = C(X)$
and $X$ is completely regular,
then $\beta(X,\SS) = \beta(X)$, and we have just given one of the
standard definition of the \v Cech compactification.
If $X$ is completely regular, then
every compactification  of $X$ is of the form $\beta(X,\SS)$ for some
$\SS$ -- namely, the collection of all those $f \in C(X)$ which 
extend to the compactification.

If $\TT \subseteq \SS \subseteq C(X)$, let us use $\pi$ to denote
the natural projection from $\beta(X,\SS)$ to $\beta(X,\TT)$.
In the case $\SS = C(X)$, this is just expressing the maximality
of $\beta(X)$ among all compactifications.  If $\TT$ 
``generates'' $\SS$, then $\pi$ is a homeomorphism.  More precisely,
let $c(\TT)$ denote the closure of $\TT$ in
the Banach algebra $C(X)$; this is the smallest closed linear
subspace of $C(X)$ containing $\TT$ and closed under pointwise
products of functions.

\medskip

{\bf Lemma 3.0.}  $\beta(X,\TT)$ and $\beta(X,c(\TT))$ are homeomorphic.

{\bf Proof.} It is easy to check that the projection
$\pi : \beta(X,c(\TT)) \to \beta(X,\TT)$ is 1-1; that is,
if $\psi,\varphi \in \beta(X, c(\TT))$ and $\psi(f) = \varphi(f)$
for all $f \in \TT$, then $\psi(f) = \varphi(f)$ for all
$f \in c(\TT)$.  \eop

\medskip

The functorial properties of these compactifications are a little
complicated because of the additional parameter, $\SS$.
Suppose that $X,Y$ are both compact spaces and $h : X \to Y$ is a continuous
function, $\SS \subseteq C(X)$, and $\TT \subseteq C(Y)$.
If we know that $h \circ f \in \SS$ for
each $f \in \TT$, then in a natural way we can define a continuous
function $\beta(h,\SS,\TT) : \beta(X,\SS) \to \beta(Y,\TT)$ by
$\beta(h,\SS,\TT)(\psi)(f) = \psi(h \circ f)$.

\smallskip

Returning now to models, 
let $M \subseteq N$ be two transitive models of $ZFC$, and
we define $\Phi=\Phi_{M,N} : CT^M \to CT^N$ as follows.
If $X \in CT^M$, let $\Phi(X)$ be $(\beta(X, C(X) \cap M))^N$.
More verbosely, working within $N$, we have the same space $X$, and
we use $C(X) \cap M $, which is a subset of $C(X)$, to compute
a compactification of $X$, which we are calling $\Phi(X)$.
This $\Phi$ is functorial in the following sense:  Let $h$
be a morphism of $CT^M$; that is, $h,X,Y \in M$ and, in $M$, $h$
is a continuous map from $X$ to $Y$, where $X,Y \in CT^M$.  Then
in $N$, $h : X \to Y$ is still continuous, and we
may extend it to $\Phi(h) : \Phi(X) \to \Phi(Y)$ by
letting $\Phi(h) = \beta(h, C(X) \cap M, C(Y) \cap M) $.
It is now easy to check from the definitions that

\medskip

{\bf Lemma 3.1.}  $\Phi_{M,N}$ is a covariant functor from
$CT^M$ to $CT^N$.

\medskip

Lemma 3.0 may be used to verify that, as claimed above,
$\Phi(X)$ is the ``obvious thing''.  The point is, we often
do not need the full $C(X) \cap M$, but may get by with some sub-class.

\medskip

{\bf Lemma 3.2.}  Suppose that in $M$, $X$ is compact,
 $\TT \subseteq C(X)$, and
$c(\TT) = C(X)$.  Then in $N$, $\Phi(X) = \beta(X, \TT)$.

{\bf Proof.}  Observe that  in $N$, 
$c(\TT) = c(C(X) \cap M)$, and apply Lemma 3.0.  \eop

\medskip

We mention two special cases of this.  First, suppose in $M$ that
$X$ is a compact subset of Euclidean space, $\RR^k$.
Let $\TT$ be the set of the  $k$ co-ordinate projections.  By
the Stone -- Weierstrass Theorem (applied in $M$), $c(\TT) = C(X)$.
But then in $N$, 
$\Phi(X) = \beta(X, \TT)$, which is just the closure of $X$ computed in
the $\RR^k$ of $N$.  In particular, if $X$ is, say, the $n$-sphere
of $M$ (so, $k = n+1$), then $\Phi(X)$ is the $n$-sphere of $N$.
Second, if $X$ is a compact zero dimensional space in $M$,
we may let $B$ be the clopen algebra of $X$, so that
$X$ is the Stone space of $B$.  In $M$, let
$\TT$ be the set of all continuous maps from $X$ into
$\{0,1\}$; then $c(\TT) = C(X)$.  From this, it is easy to see that in $N$,
$\Phi(X)$ is the Stone space of the same $B$, computed within $N$.

It is also easy to see that $\Phi$ preserves subspaces and products.
Also, if $X$ is a LOTS, then $\Phi(X)$ is the Dedekind completion
of the same LOTS; to see this, apply the above method, with $\TT$ the
set of non-decreasing real-valued functions.

See Bandlow [1] for a somewhat different treatment of $\Phi$.

\smallskip

We turn now to measures.  This is easiest to approach 
via the Riesz Representation Theorem, viewing measures as
linear functionals on $C(X)$.
If $h \in C(X)^M$, and $X \in CT^M$, then in $M$, $h$ is a continuous
map from $X$ to an interval $[a,b]$.  So, in $N$, we
we have $\Phi(h)$, which maps $\Phi(X)$ into $\Phi([a,b])$, which is
the $[a,b]$ of $N$.  So, $\Phi(h) \in C(\Phi(X))^N$.

\medskip

{\bf Lemma 3.3.}  Let $X$ be as above,
a compact Hausdorff space in $M$. In $N$, $\Phi$ is an isometric embedding
of $C(X) \cap M$ into $C(\Phi(X))$, and 
$C(\Phi(X))$ is the closed linear span of $\Phi(C(X) \cap M)$.

\medskip

In particular, suppose that in $M$, $\mu$ is a Radon measure on the
compact space $X$.  Then, via integration, $\mu$ defines a 
positive linear functional on $C(X)$, and by Lemma 3.3, this
linear functional extends uniquely to a positive linear
functional on the $C(\Phi(X))$ of $N$, which,
by the Riesz Representation Theorem,
corresponds uniquely to a Radon measure on $\Phi(X)$.
We call this measure $\Phi(\mu)$.  
Suppose, now, that in $M$, $\mu$ is non-separable.  Then, in $M$
we may find, for some fixed $\epsilon > 0$,  functions
$h_\alpha \in C(X)$ for
$\alpha < \omega_1$ such that the $L^1(\mu)$ distance between the 
$h_\alpha$ is at least $\epsilon$.  Then, this same situation will
persist in $N$ -- that is, in $N$, $L^1(\Phi(\mu))$ will be non-separable,
and hence $\Phi(\mu)$ will be non-separable, assuming that
$\omega_1$ has the same meaning in $M$ and $N$.  Thus,

\medskip

{\bf Lemma 3.4.}  Suppose that  $M \subseteq N$ are two transitive models
of $ZFC$, $\omega_1^M = \omega_1^N$,
and in $M$, $X$ is compact and $X \notin MS$.
Then in $N$, $\Phi(X) \notin MS$.

\medskip

Of course, $\omega_1^M = \omega_1^N$ is necessary.  For example, for any $X$,
if the weight of $X$ becomes countable in $N$, then $\Phi(X)$ will
be second countable in $N$ and hence be in $MS^N$.

The preservation of the property ``$X \in MS$'' is more tricky, as we
discuss in the next section.  It is quite possible that
$X \in MS^M$, but $N$ is some generic extension of $M$ which
adds a new measure which happens to be non-separable.
The forcing can even be $ccc$, in which case  $\omega_1^M = \omega_1^N$.
It is not hard to see that  ``$X \in MS$'' is preserved by any
forcing which has $\omega_1$ as a precaliber.

\smallskip

Note that for zero dimensional spaces, the results of this
section all reduce to trivialities.  If, in $M$, $X$ is the Stone space of
the Boolean algebra $B$, then $\Phi(X)$ will simply be the
Stone space of $B$ as computed in $N$.  Furthermore, if in $M$,
$\mu$ is a Radon measure on $X$, then $\mu$ is determined by its
values on the clopen sets -- i.e., by a finitely additive measure on $B$ --
and in $N$, that same finitely additive measure determines a Radon
measure, $\Phi(\mu)$,  on $\Phi(X)$.

\bigskip

{\bf \S4. Destroying Membership in $MS$.} In this section we show that
being in $MS$ can be destroyed by a $ccc$ forcing -- specifically,
by forcing with a Suslin tree.   Now, the functor $\Phi$ of the
previous section is from the $CT$ of the ground model, $V$, to
the $CT$ of a generic extension of $V$.
In the generic extension, $X$ will contain a copy of
$2^{\omega_1}$, which, by Lemma 2.0, will be sufficient to imply
that $X \notin MS$.

\medskip

{\bf Theorem 4.0.}  If there is a Suslin tree, $T$, then
there is a Corson compact space $X \in MS$ such that
$T$ forces that $X$ contains a homeomorphic copy of $2^{\omega_1}$.

{\bf Proof.}
Actually, our proof just uses the fact that $T$ is Aronszajn;
except that forcing with an Aronszajn tree does not in general
preserve $\omega_1$.  In any case, $T$ will force that
$X$ contains a homeomorphic copy of $2^{\lambda}$, where $\lambda$
is the $\omega_1$ of the ground model, but this is trivial if
$\lambda$ becomes countable in the $T$ extension.

As usual, $Lev_\alpha(T)$ denotes  level $\alpha$ of the
tree and $T_\alpha = \bigcup_{\xi < \alpha} Lev_\xi(T)$.
Let us use $\sqle$ for the tree order.

We shall construct the space $X$ from the chains of $T$.  
Identify $\PPP(T \times 2)$ with $2^{T \times 2}$ by
identifying a set with its characteristic function.  Giving
$2^{T \times 2}$ its usual topology makes
$\PPP(T \times 2)$ into a compact space.  If 
$x \in \PPP(T \times 2)$, let $\hat x \in \PPP(T)$ be its projection:
$\hat x = \{t \in T : \exists i < 2 (\langle t,i \rangle) \in x\}$.
Let $X$ be the set of all $x \in \PPP(T \times 2)$ such that
$\hat x$ is a downward-closed chain in $T$; that is
$\forall s,t \in \hat x ( s \sqle t \vee t \sqle s)$ and
$\forall t \in \hat x \forall s \sqle t (s \in \hat x)$.
Note that $X$ is closed, and hence compact.  Since $T$
is Aronszajn, each such $\hat x$ is countable; so, identifying
sets with characteristic functions, every $x \in X$ is eventually 0,
so $X$ is a compact subspace of a $\Sigma$-product of copies
of $\{0,1\}$, and hence Corson compact (see \S2).

Now, it is easy to see that in any extension, $V[G]$, of $V$,
$\Phi(X)$ is just the space defined from the same tree, by the
same definition.  However, if in $V[G]$, there is an uncountable
maximal chain $C \subseteq T$, then 
$\{x \in \PPP(T \times 2) : \hat x = C\}$ will be a subspace of
$\Phi(X)$ homeomorphic to $2^C$, which is homeomorphic to $2^{\omega_1}$.

So, we are done if we can prove (in $V$) that $X \in MS$.  
Now, one cannot prove in ZFC that every Corson compact is in $MS$ [12],
but this one is.

Let $\nu$ be a Radon measure on $X$.
For each $t \in T$, let $X_t = \{x \in X : t \in \hat x\}$.
This is closed, and hence measurable.  
If $\epsilon \ge 0$, let $T^\epsilon = \{t \in T : \nu(X_t) > \epsilon\}$.
Note that $T^\epsilon$ is a sub-tree of $T$.
If $\epsilon > 0$, then each level of $T^\epsilon$ is finite
(since the $X_t$ are disjoint for $t$ on a given level of $T$).
Since $T$ is Aronszajn, this implies that $T^\epsilon$ is countable
for each $\epsilon > 0$.  Letting $\epsilon\searrow 0$, we see
that $T^0 = \{t \in T : \nu(X_t) > 0\}$ is countable.

So, we can fix an $\alpha < \omega_1$ such that for each
$s \in Lev_\alpha(T)$,\ \  $\nu(X_s) = 0$.  Let
$F = \bigcup\{X_s : s \in Lev_\alpha(T)\}$; then $F$ is
a null set, and $X \backslash F$ is homeomorphic
to a subspace of $\PPP(T_\alpha \times 2)$, and hence second countable.
Since every finite Borel measure on a second countable space is
separable, $\nu$ is separable.   \eop

\vfill\eject

\centerline{\bf REFERENCES}

\bigskip

{
\parskip = 3pt plus1pt
\parindent = 7 true mm

\item{[1]} I. Bandlow, On the origin of new compact spaces
in forcing models, {\it Math. Nachrichten\/} 139 (1988) 185-191.

\item{[2]} M. D\v zamonja and K. Kunen, Measures on compact $HS$ spaces,
{\it Fundamenta Mathematicae\/} 143-1 (1993) pp 41-54.

\item{[3]} D. Fremlin, Consequences of Martin's Axiom, {\it Cambridge University
Press\/}, 1984.

\item{[4]} D. Fremlin, Real-valued measurable cardinals, {\it in the Israel
Math.~Conf.~Proceedings, Haim Judah (ed.)}, Vol 6. (1993), 151-304.
 
\item{[5]} R. J. Gardner and W.F. Pfeffer, Borel measures, {\it in the
Handbook of Set-Theoretic Topology, K. Kunen and
J.E. Vaughan (ed.), North-Holland\/}, 1984, 961-1044.

\item{[6]} M. Gitik and S. Shelah, Forcings with ideals and simple forcing notions,
{\it Israel J.~Math\/} 68 (1989) 129-160.

\item{[7]} M. Gitik and S. Shelah, More on simple forcing notions and forcings
with ideals, {\it Annals of Pure and Applied Logic\/} 59 (1993) 219-238.

\item{[8]} P. Halmos, Measure Theory,
{\it Van Nostrand Reinhold Company}, 1950.

\item{[9]} R. Haydon,  On dual $L^1$--spaces and injective
bidual Banach spaces, {\it Israel J.~Math\/}  31 (1978) 142--152.

\item{[10]} J. Henry, Prolongement des measures de Radon, {\it Ann.~Inst.
Fourier \/} 19 (1969) 237-247.

\item{[11]} K. Kunen,  A compact L--space under CH, {\it
Topology Appl.\/} 12 (1981) 283--287.

\item{[12]} K. Kunen and J. van Mill, Measures on Corson compact spaces,
{\it to appear}.

\item{[13]} D. Maharam,  On homogeneous measure algebras,
{\it Proc.~Nat.~Acad.~Sci.~USA\/} 28 (1942) 108--111.

\item{[14]} H. P. Rosenthal,  On injective Banach spaces and the spaces
$L^\infty (\mu)$ for finite measures $\mu$,
{\it Acta Math.\/} 124 (1970) 205--248.

}

\bye